# DESENVOLVIMENTO DE UM MODELO DE PROGRAMAÇÃO CONVEXA PARA O PROBLEMA DE FLUXO DE POTÊNCIA ÓTIMO

**Mauro Viegas da Silva[1,3], Mahdi Pourakbari-Kasmaei[2],**
**José Roberto Sanches Mantovani[3]**

[1]Universidade do Estado de Mato Grosso – UNEMAT
Departamento de Matemática
Cáceres – MT

[2]Aalto University
Department of Electrical Engineering and Automation
Espoo, Finland

[3]Universidade Estadual Paulista "Júlio de Mesquita Filho" – UNESP
Departamento de Engenharia Elétrica
Faculdade de Engenharia de Ilha Solteira
Ilha Solteira – SP

mauroviegas@unemat.br[*,1,3], Mahdi.Pourakbari@aalto.fi[2], mant@dee.feis.unesp.br[3]

**RESUMO**

O fluxo de potência ótimo (FPO) é uma classe de modelos de otimização dedicada ao desenvolvimento de ferramentas computacionais para o planejamento e operação de sistemas elétricos de potência (SEP). Neste trabalho apresenta-se um modelo convexo estendido baseado na formulação polar do FPO. Nesta formulação os termos senoidais e cossenoidais são as componentes mais difíceis do processo de convexificação do modelo matemático, porque funções desse tipo oscilam entre côncavas e convexas. Estas funções não são inicialmente convexas, mas é possível encontrar um *underestimator* convexo (UC) para as mesmas. Para obter este UC, as séries de Taylor, embora não convexas apresentam boas aproximações para essas funções trigonométricas. Embora, o modelo permaneça não convexo, esses termos podem ser reformulados para os correspondentes termos convexos equivalentes. O modelo de FPO convexo desenvolvido neste trabalho é testado e analisado nos sistemas testes IEEE de 14, 30, 54 e 118 barras.

**Palavras-chave:** Programação Convexa. Fluxo de Potência Ótimo. Séries de Taylor.

**Área principal: PO na área de Energia.**

**ABSTRACT**

The optimal power flow (OPF) is an optimization model dedicated to the development of computational tools used for the planning and operation of electric power systems (EPS). In this work, based on the polar formulation, an extended convex model is presented. To do so, the sinusoidal and cosinusoidal terms are the toughest part of the convexification process, since such types of functions oscillate between concave and convex. These functions are not initially convex, but there is a possibility of finding a convex *underestimator* (CU) for them. To obtain this CU, the Taylor series presents a good, however nonconvex, approximation for such trigonometric functions. Although the model remains nonconvex, these terms can be recast to the corresponding equivalent convex terms. The obtained convex model of the OPF is tested and analyzed using the IEEE 14-, 30-, 54-, and 118-bus test systems.

**KEYWORDS;** Convex Programming. Optimal Power Flow, Taylor series.

**Main area: PO in the energy area.**

## 1. Introdução

O fluxo de potência ótimo (FPO) é uma classe de problemas de otimização que visa obter os ajustes das variáveis de controle, para as condições de operação da rede de energia elétrica (despacho de potências ativa e reativa, redução dos custos de operação, entre outras), atendendo a um conjunto de restrições físicas e operacionais dos equipamentos e da rede de transmissão. O termo FPO, de modo geral, refere-se à otimização de um estado do sistema como o caso base (CARPENTIER, 1962) e o FPO com restrição de segurança (FPORS), que está relacionado a um processo de otimização de cenários de operação do sistema sujeitos a um conjunto adicional de restrições de segurança, ou seja, restrições de contingências, assim o FPO é um caso especial de FPORS (ALSAC e STOTT, 1974).

Nos problemas de FPO determinísticos as variáveis que representam o estado do sistema, referentes às magnitudes de tensão e ângulos nas barras, são contínuas e as variáveis que representam os "taps" de transformadores defasadores e dos transformadores com regulação de tensão, bancos de capacitores e de reatores, topologia da rede, etc são de natureza discreta, o que caracteriza o FPO como um problema real de grande porte não-linear inteiro misto (PNLIM), não convexo. Desta forma, o problema de FPO por ser um problema de otimização de grande porte e não convexo apresenta muitos desafios, devido às limitações dos métodos de programação matemática e *solvers*[1] para sua solução, bem como do limite de memória e limitação da capacidade dos processadores dos computadores.

Nos últimos dez anos o foco das pesquisas na formulação e solução do FPO tem mudado consideravelmente, visando a obtenção de modelos convexos eficientes computacionalmente e considerando a modelagem de incertezas inerentes às variáveis envolvidas no modelo de SEP reais. Assim, como uma extensão deste estudo, neste trabalho o modelo do problema de FPO tradicional de PNLIM é reformulado como um modelo *framework*[2] para otimização global que deve apresentar soluções ótimas globais e permitir que o problema seja resolvido usando *solvers* comerciais. Este modelo de programação matemática convexo contempla a redução dos custos operacionais das unidades geradoras térmicas, restrições físicas, operacionais e os ajustes dos controles de potência reativa existentes no sistema. Desta forma, propõe-se um modelo convexo para o FPO formulado na forma polar considerando apenas variáveis reais e são apresentados e discutidos os resultados de testes obtidos através da implementação computacional deste modelo no solver AMPL/KNITRO, utilizando-se os sistemas IEEE de 14, 30, 54 e 118 barras.

## 2. Revisão Bibliográfica - Relaxação Convexa FPO

A busca pela otimalidade global considerando certas hipóteses de um problema não-convexo de PNLIM como propostos por (FLOUDAS e MARANAS, 1995), (PARDALOS e ROMEIJN, 1995), (BOYD e VANDENBERGHE, 2004), (LIBERTI e MACULAN, 2006), (BAI, WEI, *et al.*, 2008) e (LAVAEI e LOW, 2012) deram início a uma nova linha de pesquisa que investiga várias relaxações convexas como em (LOW, 2014), (LOW, 2014) e (P, CAMPI, *et al.*, 2014). Quando se trabalha com relaxações convexas, estas podem mostrar informações importantes a respeito do problema de PNLIM (BOYD e VANDENBERGHE, 2004), por exemplo, o limite inferir da solução do PNLIM, o que permite verificar o grau de subotimalidade de uma solução de otimização local, mostrando que a relaxação utilizada é de alta qualidade, ou seja, *tight* (BOYD e VANDENBERGHE, 2004).

Para o *gap* de dualidade do problema relaxado igual zero, tem-se que o resultado do problema relaxado convexo é também solução do problema original, sendo o ótimo global do problema original obtido. Se o *gap* de dualidade é diferente de zero, geralmente, a solução do problema convexo relaxado não é fisicamente significativa (D, LESIEUTRE e DEMARCO, 2014). Outra característica da relaxação convexa está em certificar a factibilidade do problema. A relaxação convexa pode ser usada em diversas áreas, como por exemplo no planejamento de sistema de distribuição (ORTIZ, POURAKBARI-KASMAEI, *et al.*, 2018). Na Figura 1 ilustra-se, intuitivamente, o conceito de uma relaxação convexa.

---
[1] Um *solver* é um software matemático, sob a forma de um programa de computador autônomo ou como uma biblioteca de software, que "resolve" um problema matemático.

[2] Um *Framework* ou arcabouço conceitual é um conjunto de conceitos usado para resolver um problema de um domínio específico. *Framework* conceitual não se trata de um software executável, mas sim de um modelo de dados para um domínio. **Fonte:** https://pt.wikipedia.org/wiki/Framework**,** acessado em 02-01-2018.

*Figura 1 — Conjunto convexo relaxado a partir de um conjunto não convexo.*

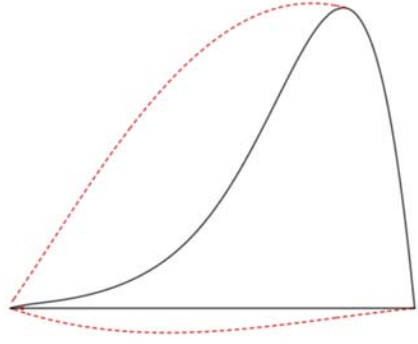

A abordagem de relaxação convexa tem apresentado teoricamente algumas vantagens em termos de confiabilidade e esforço computacional, ou seja, a solução é fornecida em um tempo polinomial.

2.1 Formulação do FPO-AC

A formulação do problema de FPO-AC é desenvolvida utilizando-se o sistema de coordenadas polares (POURAKBARI-KASMAEI, M e colab., 2016; POURAKBARI-KASMAEI, Mahdi e colab., 2014). Considere um sistema elétrico composto por $N$ barras e K ramos. Sejam $n, m$ os índices das barras e $n, m \in 1,2, \dots, N$, e seja $k$ um componente do sistema de transmissão (linhas, transformadores) com seus terminais entre as barras $n$ e $m$, $k \in 1,2, \dots, K$. O problema clássico e básico de FPO pode ser formulado como:

$$min: c(Pg_k) \tag{1}$$

$$Pg_k - Pd_k - \sum_{km} V_k V_m (G_{km} cos\theta_{km} + B_{km} sen\theta_{km}) = 0 \tag{2}$$

$$Qg_k - Qd_k - \sum_{km} V_k V_m (G_{km} sen\theta_{km} - B_{km} cos\theta_{km}) = 0 \tag{3}$$

$$V_k^{min} \leq V_k \leq V_k^{max} \tag{4}$$

$$P_{g_k}^{min} \leq P_{g_k} \leq P_{g_k}^{max} \tag{5}$$

$$Q_{g_k}^{min} \leq Q_{g_k} \leq Q_{g_k}^{max} \tag{6}$$

$$(|S_{nmk}|)^2 \leq (S_k^{max})^2 \ \forall k \tag{7}$$

Nesta formulação, as equações (2) e (3) representam $2N$ restrições de igualdades não lineares com termos quadráticos e funções envolvendo senos e cossenos. Nas simulações descritas no capítulo 5 a equação (1) como em (POURAKBARI-KASMAEI e MANTOVANI,

2018) minimizam os custos de geração de unidades térmicas, ou seja $c(Pg_k) = c_0 + c_1(Pg_k) + c_2(Pg_k)^2$, onde $Pg_k$ e $Qg_k$ correspondem às potências ativa e reativa geradas pelo gerador conectado na barra $k$, $c_0, c_1$ e $c_2$ correspondem aos coeficientes de custos. A voltagem na barra $k$ é dada por $V_k$ já as potências ativas e reativas das demandas de cargas na barra $k$ são representados por $Pd_k, Qd_k$. Seja $Y = G + jB$ a matriz complexa de admitância onde $G_{km}$ e $B_{km}$ são seus elementos. $V_k^{min}, V_k^{max}$ são as magnitudes de tensão mínima e máxima na barra $k$. $P_{g_k}^{min}, P_{g_k}^{max}$ são os limites mínimo e máximo de geração de potência ativa do gerador conectado na barra $k$. $Q_{g_k}^{min}; Q_{g_k}^{max}$ são os limites mínimo e máximo de geração de potência reativa do gerador conectado na barra $k$. $S_{nmk}$ é o fluxo de potência aparente através do componente $k$ do sistema, conectado entre as barras $n$ e $m$. $S_k^{max}$ é o fluxo máximo de potência aparente permitido no elemento $k$.

**3 Problema de Otimização**

Considere um ponto $x_0 \in X$ de uma função $f(x)$ no conjunto $X \subseteq \mathbb{R}^n$. O ponto $x_0$ é considerado um ***mínimo global*** se $f(x_0) \leq f(x)$ para todo $x \in X$. A existência de $x_0$ é garantida se $f(x)$ é contínua e se $X$ é fechado e limitado, ou seja, compacto. O ponto $x_0$ é um ***mínimo local*** de $f(x)$ se existe um $\epsilon > 0$ tal que $f(x_0) \leq f(x)$ para todo $x \in X$ com $\| x - x_0 \| \leq \epsilon$. Assim, todo mínimo de uma função convexa tem o mesmo valor da função e são, portanto, globais. De modo geral, uma função pode ter vários mínimos locais e um mínimo global.
Um problema de otimização convexa é da forma:

$$\begin{aligned} min \quad & f_0(x) \\ s.a. \quad & f_i(x) \leq 0, \quad i = 1, \dots, m \\ & a_i^T x = b_i, \quad i = 1, \dots, p \end{aligned} \quad (8)$$

onde $f_0$ e $f_1, \dots, f_m$ são funções convexas para as variáveis $x \in X$.

3.1 *Underestimator* Convexo UC

O UC é uma importante ferramenta quando se trabalha com otimização global, uma vez que a relaxação do problema original ou a aproximação UC pode ser obtida. Quando se substitui uma função não convexa em um problema por um UC pode-se usar todos os métodos disponíveis em otimização convexa para resolver a versão aproximada do problema não convexo.

**Definição 1.** A função $g(x)$ é um *underestimator* convexo de uma função $f(x)$ não convexa com $x \in C$, onde $C$ é um conjunto convexo, se:

  i.  $g(x)$ é convexa,
  ii. $g(x) \leq f(x) \, \forall x \in C$.

Para uma dada função, existem infinidades de *urderestimators* que fornecem diferentes níveis de erros de aproximações. O *underestimator* com menor erro, ou *tightest*, é chamado de envelope convexo da função (LOCATELLI, 2010). Na Figura 2 ilustra-se o conceito de relaxação para obter um UC.

***Figura 2*** — ***A curva que representa a região factível de um conjunto não convexo e sua relaxação convexa.***

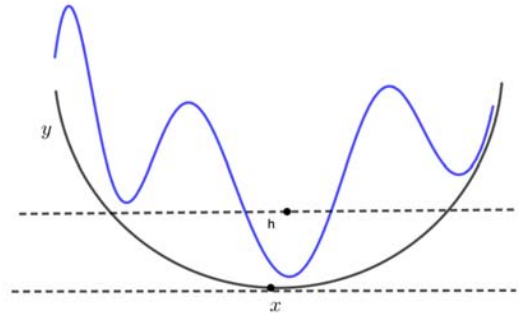

**Definição 2.** A função $g(x)$ é o *envelope convexo* de uma função não convexa $f(x)$ com $x \in C$, onde $C$ é um conjunto convexo, se

  i.  $g(x)$ é um *underestimator* convexo de $f(x)$,
  ii. $g(x) \geq h(x)\ \forall x$, para todo *underestimators* convexos $h(x)$ de $f(x)$.

Na Figura 3 estão ilustrados os conceitos desta definição.

***Figura 3*** — ***Função não convexa 1, sendo 3 e 4 underestimators convexos e 2 envelope convexo no intervalo*** $[x_1, x_2]$***.***

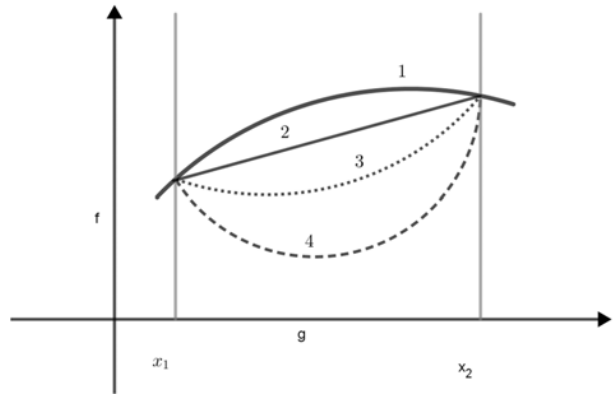

3.2 Funções Signomiais

As funções signomiais são utilizadas neste trabalho para resolver o problema de convexidade de *senos* e *cossenos* que aparecem no modelo de FPO em coordenadas polares. As signomiais constituem um grupo muito grande de funções, por exemplo, soma de polinômios e termos bi e tri-lineares são todas funções signomiais. Considere como signomiais, extensões multidimensionais de polinômios onde é permitido que as potências assumam valores não inteiros.

**Definição 3.** Uma função signomial é definida como a soma de termos signomiais, que por sua vez consiste de produtos de funções potências. Uma função signomial de $N$ variáveis e $J$ termos signomiais, pode ser expressa matematicamente como:

$$\sigma(x) = \sum_{j=1}^{J} c_j \prod_{i=1}^{N} x_i^{p_{ji}},$$

onde os coeficientes $c_j$ e as potências $p_{ji}$ são reais. As variáveis $x_i$ assumem valores reais positivos ou inteiros.

Funções signomiais são fechadas sobre adição, multiplicação escalar por constantes reais. Em geral, as funções signomiais são não convexas, porém uma condição necessária de convexidade pode ser obtida. Como a soma de termos convexos é convexa, determinando-se os termos convexos signomiais individuais pode-se obter a convexidade de toda a função.

**Teorema 1.** (Convexidade de um termo signomial positivo). Considere o termo signomial $s(x) = c \cdot x_i^{p_i} \dots x_N^{p_N}$, onde $c > 0$, o termo será convexo se uma das duas condições seguintes for satisfeita:

i. Todas as potências $p_i$ são negativas; ou,

ii. Uma potência $p_k$ é positiva, o resto das potências $p_i$, com $i \neq k$ são negativas, e

$$\sum_{i=1}^{N} p_i \geq 1,$$

ou seja, as somas das potências são maiores que ou iguais a um.

**Teorema 2.** (Convexidade de um termo signomial negativo). O termo signomial $s(x) = c \cdot x_1^{p_i} \cdots x_N^{p_N}$, onde $c < 0$, é convexo se todas as potências $p_i$ são positivas e

$$0 \leq \sum_{i=1}^{N} p_i \leq 1,$$

isto é, a soma das potências está entre zero e um.

3.5 A Reformulação *Framework* para Otimização Global

As técnicas apresentadas nesta seção encontram-se nos trabalhos de (LUNDELL, 2009) e (SKJÄL, 2014) que consistem em explorar as características das funções signomiais, ou seja, funções não convexas onde é possível obter UC para otimização global. As funções *senos* e as funções *cossenos* não são convexas e sabe-se que as séries de Taylor são ótimas aproximações para essa classe de funções. Por outro lado, série de Taylor é do tipo signomial, ou seja, é possível obter um UC, que é de grande valia para convexificar o modelo FPO desenvolvido no capítulo seguinte.

**Definição 4.** (Problema de programação signomial inteiro misto (PSIM)) Um problema de PSIM pode ser reformulado como:

$$\begin{aligned} min. \quad & f_0(\mathbf{x}) \\ s.a. \quad & A\mathbf{x} = a, \quad B\mathbf{x} \leq b, \\ & f_i(\mathbf{x}) \leq 0, \quad i = 1, \dots, m \\ & q(\mathbf{x}) + \sigma(\mathbf{x}) \leq 0, \end{aligned} \quad (9)$$

em que a função $f_0(\mathbf{x})$ é convexa, $A(\mathbf{x}) = a$ e $B(\mathbf{x}) \leq b$ são equações e inequações lineares, respectivamente. O vetor $\mathrm{x} = [x_1, x_2, \cdots, x_n]$ pode contemplar variáveis inteiras ou valores reais, ou ambos. As restrições $f_i(\mathbf{x}), i = 1, \dots, m$ são restrições de inequações não lineares que podem ser convexas ou não convexas e são adicionados $q(\mathbf{x})$ e $\sigma(\mathbf{x})$ que consistem de funções convexas e funções signomiais, respectivamente. Todas as variáveis $x_i$ nas funções signomiais $\sigma(\mathbf{x})$ são consideradas positivas.

3.3 Otimizando as Variáveis de Transformações

Um exemplo simples, para ilustrar como é aplicada a mudança de variável em um termo signomial para obter um UC unidimensional, pode ser apresentado através da função:

$$f(x) = (x^4 + 80x^2 - 160x + 100) \overbrace{-15x^3}^{(i)}, \quad 1 \leq x \leq 6,$$

de acordo com os teoremas (1) e (2) todos os termos na função $f(x)$ são convexos exceto $(i)$. Seja $x = X^p$ usando o solver *Cplex* para resolver o problema PLIM proposto (LUNDELL, WESTERLUND e WESTERLUND, 2009) obtém-se a melhor potência $p = 1/3$. Substituindo $x = X^{\frac{1}{3}}$ em $(i)$ tem-se agora uma nova função convexa.

$$f(x, X) = (x^4 + 80x^2 - 160x + 100) - 15X.$$

A potência de transformação na forma $x = X^p$ torna o termo signomial convexo para valores específicos de $p$ e transfere a não convexidade do termo signomial para uma restrição introduzida pela potência de transformação. Esta restrição está na forma de potência de transformação inversa que pode ser convexa ou não convexa, ou seja $X = x^{\frac{1}{3}}$.

Considerando o que é proposto por (LUNDELL, SKJAL e WESTERLUND, 2013) o valor de $X = x^{\frac{1}{3}}$, pode ser aproximado por uma *piece linear function* (*PLF*). Neste trabalho as novas restrições não são linearizadas, desta forma, considera-se o trabalho de (TSAI e LIN, 2008) e mais recentemente os trabalhos de (ATTARHA e AMJADY, 2016) e (AMJADY, DEHGHAN, *et al.*, 2017) que propõem um *sign transformation* para a manipulação algébrica das possíveis restrições não convexas. Na próxima seção o modelo FPO é convexificado usando séries de Taylor para obter um UC, bem como utilizam-se as técnicas de transformações para funções signomiais.

**4 Reformulação do Problema de FPO**

Uma das condições para usar as técnicas de transformações é que todas as variáveis de decisões $x = \{Pg_i, Qg_i, \theta_i, V_i\}$, $i \in N$ sejam positivas. As variáveis $Pg_i$ e $V_i$ são sempre positivas, mas $Qg_i$ e $\theta_i$ podem assumir valores negativos. Para contornar esta condição utiliza-se uma *shift transformation*, simples, mas eficiente proposta por (ATTARHA e AMJADY, 2016) que torna todas as variáveis positivas. Para tornar $Qg_i$ positiva, a seguinte manipulação algébrica é aplicada:

$$Q_{gi}^{min} = min\{Q_1^{min}, \ldots, Q_n^{min}\}, \quad n \in N \tag{10}$$

$$Q_i^{tr} = Q_{gi} + |Q_{gi}^{min}|, \quad n \in N \tag{11}$$

onde $Q_{gi}^{min}$ é o limite mínimo da potência reativa gerada na barra $i$, $Q_i^{tr}$ representa a potência reativa transformada da unidade geradora $i$. Todas as variáveis $Q_i^{tr}$ são positivas, após revolver o problema FPO, e obtido os valores $Q_i^{tr}$ as variáveis $Qg_i$ podem ser determinadas a partir da inversa de (11). Similarmente, é adicionado $2\pi$ às variáveis ângulo de fase $\theta_i$, o que faz estas variáveis tornarem-se positivas, e após a resolução do problema FPO pode-se subtrair $2\pi$ dos valores dos ângulos de fases.

O que torna o problema de FPO formulado em coordenadas polares não convexo e de difícil solução são os senos e cossenos nas equações de balanço potência, que de modo geral, caracterizam este problema multimodal, alternando sua classificação topológica entre côncava e convexa. Na Figura 4 apresenta-se um gráfico desenhado no *MatLab* de uma função $f(x_1, x_2) = sen(x_1)sen(x_2)$, com $x_1, x_2 \in [0,9]^2$. Trata-se de um exemplo simples de uma função com vários mínimos locais.

***Figura 4 — Exemplo de uma função com vários mínimos locais.***

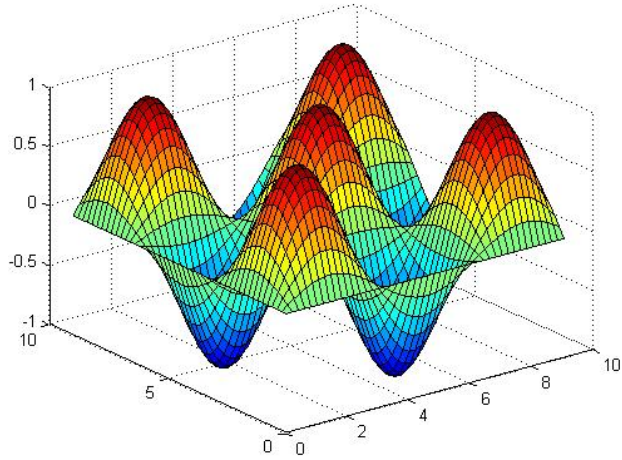

Para trabalhar com funções envolvendo *seno* e *cosseno*, tem-se da literatura (TSAI e LIN, 2008) que as séries de *Taylor* são boas aproximações, ou seja:

$$sen(x) \approx x - \frac{x^3}{3!} \qquad (12)$$

$$\cos(x) \approx 1 - \frac{x^2}{2!} \qquad (13)$$

As séries *Taylor* são funções do tipo *signomial* como apresentado na definição (3). Substituindo as funções seno e cosseno pelas respectivas expansões em séries de *Taylor* nas equações de injeção das potências ativa e reativa tem-se as novas equações (14) e (15).

$$P_n = \sum_{mk} V_n V_m \left[ G_{nmk} \left(1 - \frac{\theta_{nm}^2}{2!}\right) + B_{nmk} \left(\theta_{nm} - \frac{\theta_{nm}^3}{3!}\right) \right] \qquad (14)$$

$$Q_n = \sum_{mk} V_n V_m \left[ G_{nmk} \left(\theta_{nm} - \frac{\theta_{nm}^3}{3!}\right) - B_{nmk} \left(1 - \frac{\theta_{nm}^2}{2!}\right) \right] \qquad (15)$$

As equações de injeção de fluxos de potências (14) e (15) não envolvem *senos* e *cossenos*, porém ainda não são convexas, mas estão no formato signomial onde é possível obter um UC. Após obter este formato signomial, utiliza-se o modelo de PLIM apresentado por (LUNDELL, 2009) para gerar as potências $p_i$ que atendam as condições de convexidade, tornando as equações (14) e (15) convexas. A não convexidade é transferida para as restrições introduzidas pelas potências de transformações.

Nas restrições de balanço de potência reativa foi usado o *shift* de transformação apresentado em (10) e (11). Até esta etapa do procedimento de convexificação, o modelo ainda não está convexo, apenas foi transferida a não convexidade das equações (14) e (15) para um conjunto de restrições introduzidas pelas potências de transformações. Estas restrições estão na forma de potências de transformações inversas $X = x^{1/p}$, que por sua vez podem ser convexas ou não.

Nos trabalhos de (LUNDELL e WESTERLUND, 2009) e (LUNDELL, SKJAL e WESTERLUND, 2013) as potências de transformações inversas são aproximadas através de *PLFs*. Como as restrições geradas são monômios, utiliza-se uma proposição encontrada em (TSAI e LIN, 2008), em que uma função $f(x) = x^\alpha$, para todo $x > 0$ é convexa quando $\alpha \leq 0$ ou $\alpha \geq 1$. E $f(x)$ é côncava quando $0 < \alpha < 1$. Nos trabalhos de (AMJADY, DEHGHAN, *et al.*, 2017) e (ATTARHA e AMJADY, 2016) propõe-se uma transformação de sinal como especificada em (16).

$$\begin{cases} X_i = x_i^{\frac{1}{Q_i}}, & se \quad \frac{1}{Q_i} < 0 \quad ou \quad \frac{1}{Q_i} \geq 1, \\ X_i = -x_i^{\frac{1}{Q_i}}, & se \quad \quad 0 < \frac{1}{Q_i} < 1, \end{cases} \quad (16)$$

Assim, se uma potência de transformação inversa não é convexa ela é totalmente côncava, e pode ser facilmente convexificada mudando seu sinal.

**5 Resultados**

Os testes foram realizados usando os *softwares* Matpower 5.1 de (ZIMMERMAN, SANCHEZ e THOMAS, 2011) e (MURILLO-SÁNCHEZ, ZIMMERMAN, *et al.*, 2013) com o Matlab e Knitro com o AMPL (FOURER, GAY e KERNIGHAN, 2003). Os testes foram executados em um computador Dell Dell PowerEdge t430, sistema operacional: Linux, 2x Processadores Intel E5-2650 v4 2.2GHz, 30M Cache, 9.60 GT/s Q P I, Turbo, HT. 12Cores/24Threads (105W) Max M em 2400Hz.Todos os testes foram realizados para os seguintes sistemas testes do IEEE de 14, 30, 57 e 118 barras.

Os resultados obtidos com o modelo proposto, comparados com os resultados obtidos através do Matpower estão apresentados nas Tabelas 1, 2, 3 e 4.

*Tabela 1 — Teste no sistema IEEE 14 barras*

|  | **Modelo Não Linear** | **Modelo Convexo** |
|---|---|---|
| *Software Utilizado* | MATPOWER | KNITRO |
| Função Objetivo ($/h) | 8.081,53 | 8.081,14 |
| Tempo CPU (s) | 1,79 | 0,281 |
| Potência Ativa Gerada (MW) | 268,29 | 268,28 |
| Potência Reativa Gerada (MVAr) | 67,630 | 63,734 |
| Carga (MW) | 259,00 | 259,00 |

Fonte: Autores

*Tabela 2 — Teste no sistema IEEE 30 barras.*

|  | **Modelo Não Linear** | **Modelo Convexo** |
|---|---|---|
| *Software Utilizado* | MATPOWER | KNITRO |
| Função Objetivo ($/h) | 8.906,14 | 8906,13 |
| Tempo CPU (s) | 1,93 | 2,168 |
| Potência Ativa Gerada (MW) | 295,15 | 295,14 |
| Potência Reativa Gerada (MVAr) | 113,94 | 113,93 |
| Carga (MW) | 283,40 | 283,40 |

Fonte: Autores

*Tabela 3 — Teste no sistema IEEE 57 barras.*

|  | **Modelo Não Linear** | **Modelo Convexo** |
|---|---|---|
| *Software Utilizado* | MATPOWER | KNITRO |
| Função Objetivo ($/h) | 41.737,79 | 41.737,68 |
| Tempo CPU (s) | 0,08 | 14 |
| Potência Ativa Gerada (MW) | 1267,31 | 1267,26 |
| Potência Reativa Gerada (MVAr) | 270,56 | 270,52 |
| Carga (MW) | 1250,80 | 1250,80 |

Fonte: Autores

*Tabela 4 — Teste no sistema IEEE 118 barras.*

|  | **Modelo Não Linear** | **Modelo Convexo** |
|---|---|---|
| *Software Utilizado* | MATPOWER | KNITRO |
| Função Objetivo ($/h) | 129.660,70 | 129.660,60 |
| Tempo CPU (s) | 1,59 | 15 |
| Potência Ativa Gerada (MW) | 4.319,40 | 4319,35 |
| Potência Reativa Gerada (MVAr) | 388.26 | 388,25 |
| Carga (MW) | 4242,00 | 4242,00 |

Fonte: Autores

Analisando-se os resultados das Tabelas 1-4 verifica-se a eficiência do modelo convexo proposto para o problema de FPO em que este modelo resolvido através do KNITRO obteve as mesmas soluções obtidas pelo modelo de programação não linear resolvido através do MATPOWER. O KNITRO é desenvolvido para encontrar soluções de otimalidade local de problemas de otimizações contínuos. Uma solução local é um ponto factível, no qual o valor da função objetivo naquele ponto seja tão bom ou melhor que em qualquer ponto factível "próximo". Uma solução de otimalidade global é uma solução que fornece o melhor (isto é, o menor se o problema for de minimizar) valor da função objetivo entre todos os pontos factíveis. Se o problema é ***convexo*** todas as soluções ótimas locais são também soluções de otimalidade global.

O modelo convexo possui um número maior de restrições em função das manipulações algébricas necessárias para sua obtenção, inclusive restrições de igualdade, e isto torna maior o tempo de CPU necessário para o solver obter a sua solução ótima. Apesar de o modelo convexo consumir maior tempo computacional devido estas manipulações utilizadas para convexificar o problema original, ele garante a solução ótima global e permite concluir sobre a qualidade das soluções obtidas pelo modelo não convexo. Nos resultados obtidos até o presente desenvolvimento deste trabalho verifica-se que tanto o modelo não convexo quanto ao convexo, fornecem as mesmas soluções ótimas para todos os sistemas testes simulados.

**6. Conclusões**

O problema de FPO é um problema não linear e não convexo classificado como *NP-hard* o que torna difícil a sua solução são as equações de balanço de potência ativa e reativa que são restrições de igualdade e que apresentam algumas características especiais, tais como, a presença de senos e cossenos em sua estrutura. O modelo mais completo de FPO envolvendo variáveis binárias e discretas bem como considerar as curvas de capabilidade térmica dos diferentes tipos de geradores pode se apresentar ainda mais complexo. Nos últimos 10 anos são encontrados diversos trabalhos na literatura para formular e resolver o problema de FPO como um modelo convexo relaxado, por programação SDP, programação quadrática que incluem as cônicas de segunda ordem, modelos linearizados, bem como o modelo *framework*.

As pesquisas sobre relaxação convexa aplicada ao problema de FPO ainda são recentes, desta forma, ainda não está clara qual a melhor relaxação convexa para este tipo de problema. Pode-se concluir que esse novo ramo de pesquisa ainda propicia o desenvolvimento de diferentes linhas de pesquisa tanto em termos de desenvolvimento de modelos como de técnicas de soluções. O tipo de relaxação proposta neste trabalho se justifica, pois ainda foi pouco explorada no problema de FPO, e o seu desenvolvimento como um problema de FPO convexo constitui uma contribuição importante na área de otimização de sistemas, formulação e solução do problema de FPO.

**Agradecimento**